\documentclass[11pt]{article}
\usepackage{amsmath,amsthm,amsfonts,amssymb,mathrsfs,epsfig,bm}
\usepackage[usenames]{color}
\usepackage[colorlinks=true]{hyperref}

\parskip        \smallskipamount

\theoremstyle{remark}

\def\Z{\mathbb{Z}}

\def\R{\mathbb{R}}

\def\P{\mathbb{P}}
\def\E{\mathbb{E}}

\renewcommand{\phi}{\varphi}
\renewcommand{\epsilon}{\varepsilon}

\newcommand{\1}{{\text{\Large $\mathfrak 1$}}}

\renewcommand{\d}{\text{\rm\,d}}
\newcommand{\eqdist}{\stackrel{\text{(d)}}{=}}

\def\vol{\operatorname{vol}}

\definecolor{mygray}{gray}{0.9}

\definecolor{deeppink}{RGB}{255,20,147}

\long\def\symbolfootnote[#1]#2{\begingroup
\def\thefootnote{\fnsymbol{footnote}}\footnote[#1]{#2}\endgroup}
\newcommand{\keywords}[1]{ \noindent {\footnotesize
             {\small \em Keywords and phrases.} {\sc #1} } }
\newcommand{\ams}[2]{  \noindent {\footnotesize
             {\small \em AMS {\rm 2000} subject classifications.
             {\rm Primary {\sc #1}; secondary {\sc #2}} } } }
\def\P{\mathbb P}
\def\E{\mathbb E}
\def\B{\mathbb B}
\def\d{\mathrm{d}}
\def\CAPRAND{\operatorname{\sf T}}
\def\erf{\operatorname{erf}}
\def\erfc{\operatorname{erfc}}
\def\realpart{\operatorname{Re}}

\def\d{\text{d}}

\begin{document}
\title{\bf The mobile Boolean model: an overview and further results
}
\author{\sc Takis Konstantopoulos}
\date{\small \em 11 July 2013 \rm [sic\footnote{A version
of a paper which is scheduled to finish in the future.}]}
\maketitle

\begin{abstract}
This paper offers an overview of the mobile Boolean stochastic geometric
model which is a time-dependent version of the ordinary Boolean model
in a Euclidean space of dimension $d$.
The main question asked is that of obtaining the law of the
detection time of a fixed set. We give various ways of thinking
about this which result into some general formulas. The formulas
are solvable in some special cases, such the inertial and Brownian 
mobile Boolean models. In the latter case, we obtain some expressions
for the distribution of the detection time of a ball, when the
dimension $d$ is odd and asymptotics when $d$ is even.
Finally, we pose some questions for future research.

\vspace*{2mm}
\keywords{Boolean model, stochastic geometry, capacity,
Brownian motion, hitting times, heat equation, Bessel process,
modified Bessel functions.}

\vspace*{2mm}
\ams{60D05, 60G55, 60J65}{60K40, 33C10}
\end{abstract}

\section{Introduction}
Particles start from the points of a point process in $\R^d$ and perform
i.i.d.\ stochastic motions. Each particle carries a set (detection set) along
(typically, a ball of fixed radius). The goal is to find the time
until one of the particles detects a fixed set $K$.
This is a paper that should have been written long ago. Its first version,
\cite{KKP03} appeared in the proceedings of an obscure conference without many
proofs. After actually working out the Brownian case, we discovered
that Spitzer \cite{SPI} had done most of the job for $d=2$ or $3$. 
A few years after that, we posed some open questions 
\cite{K4B}, some of which were resolved
by Peres {\em et al.} \cite{PSSS}.
The purpose of this paper is twofold.
First, because we think that an overview of the topic is needed.
Second, because there are several interesting questions that can be asked.
This overview contains some new elements too. For example, 
it contains a method for computing the exact distribution of the
detection time in odd dimensions and asymptotics of the
distribution in even dimensions.
We also compare what happens when we run the particles a Brownian motions
vs.\ linear motions with random speeds. 
Roughly speaking, Brownian motions discover small objects
faster than linear motions.

The paper is organized as follows. First, we give an overview of
the general Boolean model. Then we pass on to the dynamic case. 
We then pose the detection problem and provide some general formulas
involving capacity notions. The inertial Boolean model (particles move on 
straight lines with constant speed) is presented next.
The Brownian Boolean model is then analyzed in detail.
We conclude with some open problems.
The first part of the paper, i.e., up to and including Section 4, is 
largely an overview but with organized notation chosen so that the
model can be explained in much greater generality than the ``solvable'' one
(the Brownian Boolean model)
Section 5 is new (inertial Boolean model), but quite simple. 
Section 6 deals with the Brownian Boolean model and explains the new
formulas regarding the distributions of detection times and their expectations,
separately in even and odd dimensions.

Throughout the paper we use the following notations: 
$f(x) \sim g(x)$, as $x \to \infty$, means
$f(x)/g(x) \to 1$, as $x \to \infty$. Similarly, when $x$ tends to
another value. Also, $f(x) 
\stackrel{\log}{\text{\Large $\sim$}}  g(x)$
means $\log f(x) \sim \log g(x)$.
The closed unit ball in $\R^d$ centered at the origin is denoted by $\B$.
The Euclidean norm of $x \in \R^d$ is denoted by $|x|$ and the $d$-dimensional
Lebesgue measure of a Borel set $A \subset \R^d$ is denoted as $\vol(A)$.
In particular, we let $\omega_d := \vol(\B)$, keeping in mind that
\[
\omega_d = \frac{\pi^{d/2}}{\Gamma(\frac{d}{2}+1)}, \quad d \in \Z_+,
\]
which gives
\[
\omega_{2n} = \frac{\pi^n}{n!}, \quad
\omega_{2n+1}= \frac{2^{n+1} \pi^n}{(2n+1)!!}, \quad
n \in \Z_+,
\]
where  $(2n+1)!! = 1 \cdot 3 \cdot 5 \cdots (2n+1)$.
When $A, B$ are sets (perhaps random) in $\R^d$, then $A\pm B$ is their Minkowski 
sum (difference). 
For example, $A$ could be a singleton and $B$ the trajectory of a stochastic
process. The modified Bessel function of the second kind and of order $\nu \in \R$ 
is denoted by $K_\nu$. It admits the integral representation
    \[
K_\nu(z) = \int_0^\infty \exp(-z\cosh t) \cosh(\nu t) \,\d t,
\quad \realpart(z) > 0.
\]
Alternative formulas for $K_\nu$ when $\nu$ is a half-integer
are in the appendix.
We remark also that the term ``Boolean model'' is not standard. Some people
refer to it as germ-grain model and allow the underlying point process to
be general. In this paper, we reserve the term ``Boolean model'' for the
case when the underlying point process is (homogeneous) Poisson 
and use ``general Boolean model'' for the case where the point process is more
general.

\section{A little background on the general Boolean model}
The general Boolean model, also known as the germ-grain model 
\cite[Ch.\ 4]{SW08},
 is one the basic objects of study of stochastic geometry.
Consider a random point process $\Phi$ on $\R^d$ (the ``{\em germs}'')
and a random compact set $G$ (a ``{\em grain}'').
Conditionally on $\Phi$, let $\{G_x, x \in \Phi\}$ be i.i.d.\ 
copies\footnote{We think of $\Phi$ both as a random discrete set,
and as a random point process. Thus $\{x\in \R^d: x \in \Phi\}$ stands for
the random set or for the support of the random measure. When $B$ is a Borel
subset of $\R^d$ we use $\Phi(B)$ to denote the value of the
random measure at $B$, i.e., the number of points of the random set
in $B$.} of $G$.
The {\em general Boolean model} is then defined 
by\footnote{When $A, B \subset \R^d$, we let
$A+B := \{a+b:a\in A, b \in B\}$ (Minkowski addition). 
Similarly, $A-B:=\{a-b:a\in A, b \in B\}$.
We let $x+B := \{x\}+B$.}
\begin{equation}
\label{booleanmodel}
\Xi:=\bigcup_{x \in \Phi} (x+G_x).
\end{equation}
(If $G$ is a deterministic set, then $\Xi$ can also be expressed as
$\Xi = \Phi+G$.)
Alternatively, if $\mathcal K$ is the collection of compact subsets of $\R^2$,
we may consider the point process 
\begin{equation}
\label{mppN}
\mathcal N:= \{(x, G_x): x \in \Phi\}
\end{equation}
as a discrete random subset of the product space $\R^d \times \mathcal K$.
In fact, it is also a marked point process because to each $x$ there is
a unique $G_x \in \mathcal K$ such that $(x, G_x)$ is an element of the random
set $\mathcal N$. Putting an appropriate probability measure on the set of
marked point processes on $\R^d \times \mathcal K$ gives another
way of constructing a general Boolean model.

The {\em capacity functional} of the general Boolean model is defined as
\[
\CAPRAND_{\Xi}(K) := \P(K \cap \Xi \neq \varnothing),
\]
for $K$ a compact set.

The capacity functional
 is a fundamental object in the theory of random
sets \cite{MOL10}. 
If $X$ is a random locally compact subset of $\R^d$ then 
\[
\CAPRAND_X(K) := \P(K \cap X \not = \varnothing)
\]
is defined on compact sets $K$.  It is submodular, i.e.,
$\CAPRAND_X(K_1 \cup K_2) \le \CAPRAND_X(K_1) + \CAPRAND_X(K_1) - \CAPRAND_X(K_1 \cap K_2)$
and upper semicontinuous, i.e., $\CAPRAND_X(K_n) \downarrow \CAPRAND_X(K)$ 
whenever $K_n$ is a decreasing sequence of compact sets with
$\bigcap_n K_n = K$.
An example of a capacity functional is the one given above. Another 
example is obtained by considering a Brownian motion $\xi:=\{\xi(t),
t \ge 0\}$ in $\R^d$ started from some $\xi(0) \in \R^d$. 
Letting 
\[
\xi^{(t)} := \big\{\xi(s),~ 0 \le s \le t\big\}
\]
be the initial segment of $\xi$ up to time $t$, we have
$\CAPRAND_{\xi^{(t)}}(K) = \P(\xi^{(t)} \cap K \neq \varnothing)$.
If we set
\[
T_K := \inf\{t \ge 0: \xi(t) \in K\},
\]
we obtain
\[
\CAPRAND_{\xi^{(t)}}(K) = \P(T_K \le t).
\]
Thus, $t \mapsto \CAPRAND_{\xi^{(t)}}(K)$
is the distribution function of the random variable $T_K$.

Going back to the general Boolean model, we observe that
if (the law of) $\Phi$ is invariant under translations and ergodic, and
if $\vol$ denotes the Lebesgue measure on $\R^d$, then 
$C_{\Xi}(\{0\})$  is the a.s.\ (and in $L^1$) limit of 
$\vol(\Xi \cap [-h,h]^d)/h^d$, as $h \uparrow \infty$ and is
known as the {\em volume fraction} of $\Xi$.

An important particular case is when $\Phi$ is a {\em spatially homogeneous
Poisson process} with intensity $\lambda$. In this case, $\Xi$ is referred to as
a Boolean model\footnote{The terminology 
is highly nonstandard. Some authors use the term ``Boolean model''
or ``germ-grain model'' for what we called ``general Boolean model''.
We decided to reserve the term ``Boolean model'' for
the case where $\Phi$ is a Poisson process.} 
\cite{SKM08}.
The Boolean model has several computational advantages.
For example, owing to the thinning property of a Poisson process, 
the set of points $x \in \Phi$ such that $x+G_x$ intersects a given
closed set $A$,
\[
\Phi_A:=\{x \in \Phi:~ A\cap (x+G_x) \neq \varnothing\},
\]
forms an inhomogeneous Poisson process in $\R^d$ with 
intensity measure
\[
\E \Phi_A(\d x) =
\lambda_A(x)\, \d x = \lambda \, \P(A \cap (x+G) \neq \varnothing)\,\d x.
\]
Indeed, conditionally on $\Phi$, the random
variables $\big\{ \1(A \cap (x+G_x) \neq \varnothing),\, x \in \Phi\big\}$
are independent and $\1(A \cap (x+G_x) \neq \varnothing)
\eqdist \1(A \cap (x+G) \neq \varnothing)$, for all $x \in \Phi$.
Since $A\cap (x+G) \neq \varnothing$ is equivalent to 
$(A-x) \cap G \neq \varnothing$, we have
\[
\lambda_A(x) = \lambda \CAPRAND_G(A-x).
\]
Notice also that if $\Xi$ is a Boolean model then the marked
point process $\mathcal N$ introduced in \eqref{mppN} is a Poisson
process on $\R^d \times \mathcal K$ with intensity measure
$\lambda\, \d x\, \P(G \in \d g)$ on the space 
$\R^d \times \mathcal K$.
Conversely, given a finite measure $\mu(\d g)$ on $\mathcal K$, we can
construct a Poisson process $\mathcal N$ on 
$\R^d \times \mathcal K$ with intensity measure $\d x\, \mu(\d g)$. 
The projection of $\mathcal N$ on $\R^d$ is the set of germs and
the projection on $\mathcal K$ is the set of grains.

In applications of stochastic geometry, it is sometimes the case that the grains
depend on a parameter $t$ (time) and increase with $t$. 
The general mobile Boolean 
model has a time parameter $t$ which principally affects the locations 
of the germs.

\section{The general mobile Boolean model}
\subsection{Definitions}
In its most general form, the general mobile Boolean model is defined
in terms of a germ point process $\Phi$ on $\R^d$ and of 
a time-dependent grain process $G=\big\{ G(t), t \ge 0 \big\}$.
The latter is assumed to be a stochastic process with values in
$\mathcal K$ and continuous sample paths, when $\mathcal K$ is 
equipped with the topology induced by the Hausdorff metric.
Conditionally on $\Phi$, let $\big\{ G_x, x\in \Phi \big\}$ be i.i.d.\ copies
of $G$ and define, for each $t \ge 0$, the general Boolean model
\[
\Xi(t) = \bigcup_{x \in \Phi} (x+G_x(t)).
\]
The {\em general mobile Boolean model} is the 
{\em process} $\Xi=\big\{
\Xi(t), t \ge 0 \big\}$.
Physically, we think of particles located at points $x \in \Phi$ at time
$0$, each having a grain $G_x(0)$ ``around it''. At time $t$, point $x$
moves to a new position and, as a result, the grain around it becomes
$x+G_x(t)$.
Notice that the set
\[
\mathcal W(t) := \bigcup_{0 \le s \le t} \Xi(s)
\]
represents everything covered by the general Boolean model up to time $t$
and is itself too
a general Boolean model with grains the points of $\Phi$
and germs being i.i.d.\ copies of 
\[
G^{(t)}:=\bigcup_{0 \le s \le t} G(s).
\]
Indeed,
\begin{equation}
\label{Wexpression}
\mathcal W(t) = \bigcup_{x \in \Phi} (x + \bigcup_{0 \le s \le t} 
G_x(s)) = \bigcup_{x \in \Phi} (x + G_x^{(t)}).
\end{equation}
Note that $\mathcal W(t)$ is increasing  in $t$.

\subsection{The canonical form}
\label{cfsection}
The canonical way of introducing the distribution of the grain process 
$\big\{ G(t), t \ge 0\big\}$
is by considering a specific point $h(G(0))$ of $G(0)$
as its ``center'' (where, formally, $h$ is a measurable function from
$\mathcal K$ into $\R^d$). 
Then
\[
\xi(t) := h(G(t)), \quad t \ge 0,
\]
describes the motion of the center.
We may thus describe the law of $\big\{ G(t), t \ge 0\big\}$, in two steps.
First by specifying the law of $\big\{ \xi(t), t \ge 0 \big\}$ and then, 
conditionally on $\big\{ \xi(t), t \ge 0 \big\}$, by specifying the law
of 
\[
D(t):=G(t)-\xi(t),\quad t \ge 0.
\]
(In the simplest case, we may assume that $\big\{ D(t), t \ge 0\big\}$
is independent of $\big\{ \xi(t), t \ge 0 \}$.)
We can then write 
\[
\Xi(t) = \bigcup_{x \in \Phi} (x+\xi_x(t) + D_x(t)),
\]
where $\xi_x(t) = h(G_x(t))$, $D_x(t)=G_x(t)-\xi_x(t)$.
The two representations for $\Xi(t)$ are absolutely equivalent. Notice,
however, that  if we let
\[
\Phi(t) :=\big \{x + \xi(t):~ x \in \Phi \big \}
\]
we can write, with a slight abuse of notation as regards the
indices of $D_x(t)$,
\[
\Xi(t) = \bigcup_{x \in \Phi(t)} (x + D_x(t)),
\]
and think of $\Xi(t)$ as a general Boolean model with
germs the points of $x\in \Phi(t)$ and grains $D_x(t)$.

In the simplest case, $\Phi$ is a homogeneous Poisson process
with intensity $\lambda$,
and $\{D(t), t \ge 0\}$
is independent of $\{\xi(t), t \ge 0\}$.
Then $\Phi(t)$ is again a homogeneous Poisson process 
with intensity $\lambda$.
Thus, if $\Xi(0)$ is a Boolean model and the trajectory of the center
is chosen independently of $D$, then,
for each $t>0$, 
$\Xi(t)$ is also a Boolean model identical in distribution to $\Xi(0)$.

Note: Without further ado we shall, henceforth, define a mobile
Boolean model to be one for which the germ point process $\Phi$ is
homogeneous Poisson and the independence between $\xi(\cdot) = h(G(\cdot))$
and $D(\cdot)$ holds.
Dropping the adjective homogeneous gives an inhomogeneous mobile Boolean model.
An inhomogeneous mobile Boolean model is arguably a good model for
a sensor network. Sensing devices are initially located at the points $x$ of
the inhomogeneous Poisson process $\Phi$ and move according to independent
random motions $x+\xi_x(t)$ (assuming $\xi_x(0)=0$). The set $D_x(t)$
represents the part of space which can be sensed by the sensor at time $t$.
If the problem is to discover an unknown target, then inhomogeneity in $\Phi$
allows for the possibility of incorporating prior information about
the location of the target. Randomness in $D_x(t)$ may model the 
different sensing abilities of the devices. And time-dependence
in $D_x(t)$ allows for modeling of loss  (or gain) of energy of the device.
Finally, randomness in the trajectories is natural too.

\section{The detection problem}
We now consider the problem of finding the distribution of
the first time that a general mobile Boolean model will detect a fixed
compact set $K$. The expressions \eqref{detGENERAL} and \eqref{detINHOMmbm}
below concern, respectively, a general mobile Boolean model and an inhomogeneous
mobile Boolean model.
The expressions \eqref{detHOMmbm}, and \eqref{detHOMmbmCANON}
concern both a homogeneous mobile Boolean model. The last
one relates the distribution of the detection time of $K$
by to the distribution of the first hitting time of a set
by the process $x+\xi(t)$ representing the
random motion of a sensing device initially located at the point $x$.

\subsection{Detection time for a general mobile Boolean model}
Let $K$ be a compact subset of $\R^d$ and let
\[
S_K := \inf\{t \ge 0:~ K \cap \mathcal W(t) \neq \varnothing\}.
\]
This is called the {\em detection time of the set $K$}.
We are interested in deriving information about the law of $K$.
Under natural assumptions on the law of $\Phi$ (e.g., if $\Phi$
is a Poisson process), and
if the grain $G$ has nonempty interior with positive probability,
then, due to compactness, the probability that $K$ is contained in
$\mathcal W(0)$ is positive.
Hence $\P(S_K=0)$ is typically (e.g., under the previous assumptions) positive.
By the monotonicity of $\mathcal W(t)$, 
\begin{equation}
\label{detGENERAL}
\P(S_K \le t) = \P(K \cap \mathcal W(t) \neq \varnothing)
= \CAPRAND_{\mathcal W(t)}(K),
\end{equation}
and, by the expression \eqref{Wexpression} for $\mathcal W(t)$,
\begin{align*}
\CAPRAND_{\mathcal W(t)}(K)
&= \P\big(\bigcup_{x \in \Phi} K \cap (x+ G^{(t)}_x) \neq \varnothing\big)
\\
&= \P(\exists x \in \Phi  ~ K \cap (x+ G^{(t)}_x) \neq \varnothing).
\end{align*}
So, if we consider the point process
\begin{equation}
\label{thinnedpoisson}
\Phi_K:=\{x \in \Phi: ~ K \cap (x+ G^{(t)}_x) \neq \varnothing)\},
\end{equation}
we have
\[
\CAPRAND_{\mathcal W(t)}(K) = \P(\Phi_K \neq 0),
\]
which is the probability that  for some point $x$ of $\Phi$
the set $x+G_x^{(t)}$ intersects $K$.
This point process depends
on the general Boolean model in a rather complicated way.

\subsection{Detection time for a (possibly inhomogeneous) mobile
\label{dtimbm}
Boolean model}
Things become simple in the case of an inhomogeneous mobile Boolean model
where $\Phi$ is a Poisson process with intensity measure $\lambda(\d x)$.
In this case, arguing as earlier, $\Phi_K$ is also a Poisson process
with intensity measure
\[
\P(G^{(t)} \cap(K-x) \neq \varnothing)\, \lambda(\d x) 
=\CAPRAND_{G^{(t)}}(K-x)\, \lambda(\d x),
\]
and thus,
\begin{equation}
\label{detINHOMmbm}
\P(S_K \le t )= \CAPRAND_{\mathcal W(t)}(K) = 1-\exp
\bigg(-\int_{\R^d} \CAPRAND_{G^{(t)}}(K-x)\, \lambda(\d x)\bigg).
\end{equation}
Since $(G^{(t)} + D) \cap (K-x) \neq \varnothing$
iff $x \in K-D-\xi(s)$ for some $s \le t$, we have,
by Fubini's theorem,
\begin{equation}
\label{byfubini}
\P(S_K > t) 
=
\exp\bigg(-\E\, \lambda\bigg( \bigcup_{0 \le s \le t} [K-D-\xi(s)]\bigg)\bigg)
\end{equation}
If $\Phi$ is a spatially homogeneous Poisson process with
\[
\lambda(\d x) = \lambda \cdot \d x,
\]
then
\begin{equation}
\label{detHOMmbm}
\P(S_K > t) = 
\exp\bigg(- \lambda\, \E\vol\bigg( \bigcup_{0 \le s \le t} [\xi(s) + D-K] 
\bigg)\bigg).
\end{equation}

We point out that the sets involved in this union 
are formed by translating the set $D-K =\{x-y: x\in D, y \in K\}$
by vectors $\xi(s)$.
Put it otherwise, a particle performing motion $\xi$ carries
a neighborhood with shape $D-K$. The set swept by
the particle up to time $t$ can be called the ($D-K$)-{\em sausage} of
$\xi$ up to time $t$. The term {\em Wiener sausage} is reserved for the
case when $\xi$ is a Brownian motion.

\subsection{Detection time for a homogeneous
mobile Boolean model in its canonical form}
Suppose now that
\[
G(t) = \xi(t) + D,
\]
where $D$ is a fixed deterministic compact set, for instance a closed ball,
and $\{\xi(t),\, t \ge 0\}$
a random process with continuous sample paths and $\xi(0)=0$;
see Section \ref{cfsection}.
Then 
\[
G^{(t)} = \bigcup_{0 \le s \le t} (\xi(s) + D) = \xi^{(t)}+D,
\]
and
\[
G^{(t)} \cap (K-x) \neq \varnothing 
\iff
(\xi^{(t)}+D) \cap (K-x) \neq \varnothing 
\iff
(x+\xi^{(t)}) \cap(K-D) \neq \varnothing.
\]
Therefore, if we define the first hitting time 
\[
T^x_B := \inf\{t \ge 0:\, x+ \xi(s) \in B\}
\]
of a closed set $B$ by the process $x+\xi(\cdot)$, we have
\[
\CAPRAND_{G^{(t)}}(K-x) = \P(T^x_{K-D} \le t).
\]
Using \eqref{detINHOMmbm} and \eqref{detHOMmbm} we  obtain
\begin{equation}
\label{detHOMmbmCANON}
\P(S_K > t) = \exp\bigg(-\lambda \int_{\R^d} \P(T^x_{K-D} \le t) \,\d x\bigg).
\end{equation}
In particular, if $T^x_{K-D}$ has density $f^x_{K-D}(t)$ then
the hazard rate $h_K(t)$ of $S_K$, defined by
\[
\P(S_K \le t+\delta \mid S_K > t) = h_K(t)\delta+o(\delta),
\quad\text{ as } \delta \downarrow 0,
\]
exists and  is given by
\[
h_K(t) = \lambda \int_{\R^d} f^x_{K-D}(t)\, \d x.
\]
Computing the distribution of $S_K$ exactly may be hard, but 
asymptotics may be possible, via knowledge of the Laplace transform
of $T^x_{K-D}$ and Tauberian theorems.

\subsection{The isotropic case}
Suppose that the process $\{\xi(t),\, t \ge 0\}$, with $\xi(0)=0$, is
isotropic, i.e., that if $Q$ is a proper rotation of $\R^d$ then
$\{Q\xi(t),\, t \ge 0\}$ has the same law as $\{\xi(t),\, t \ge 0\}$.
Assuming further that
\begin{equation}
\label{Kset}
D =  r \mathbb B , \quad K = r_0 \mathbb B,
\end{equation}
where $\mathbb B = \{x\in\R^d:\, |x|\le 1\}$ is the unit ball
of radius 1 centered at the origin,
then, clearly, the integral in \eqref{detHOMmbmCANON} can be
simplified. Indeed, with $e_1=(1,0,\ldots,0)$,
\[
T^x_{K-D} = \inf\{t \ge 0:\, |x+\xi(t)| \le r+r_0\}
\eqdist \inf\{t \ge 0:\, \big|\,|x| e_1+\xi(t)\,\big| \le r+r_0\},
\]
so, letting
\begin{equation}
\label{Trrr}
T^{\rho}_{r+r_0} := \inf\{t \ge 0:\, |\rho e_1 +\xi(t)| \le r+r_0\}
\end{equation}
be the first hitting time of $r+r_0$ by the radial process
$|\rho e_1 +\xi(t)|$, we have
\begin{align*}
\int_{\R^d} \P(T^x_{K-D} \le t) \,\d x
&= \sigma_{d-1} \int_0^\infty \P(T^{\rho}_{r+r_0} \le t)\, \rho^{d-1} \,\d\rho
\\
&= \omega_d (r+r_0)^d + \sigma_{d-1}\int_{r+r_0}^\infty 
\P(T^{\rho}_{r+r_0} \le t)\, \rho^{d-1} \,\d\rho,
\end{align*}
where $\omega_d = \pi^{d/2}/\Gamma(1+d/2)$ 
is the $d$-dimensional Lebesgue measure of
the unit ball in $\R^d$ and $\sigma_{d-1} = d\omega_d$.
Hence
\begin{equation}
\label{SKrot}
\P(S_K > t) 
= e^{-\lambda \omega_d (r+r_0)^d}
\exp \bigg(-\lambda\, \sigma_{d-1} \int_{r+r_0}^\infty 
\rho^{d-1} \, \P(T^\rho_{r+r_0} \le t) \,\d \rho \bigg).
\end{equation}
Although the integral has been reduced from a $d$-dimensional one to
$1$-dimensional, finding the distribution of \eqref{Trrr} is
still a $d$-dimensional problem.

\section{The inertial Boolean model}
\label{inbm}
We take $\xi$ to be a linear stochastic motion. Let $v$ be
a random variable in $\R^d$ and let 
the motion of a typical particle be
\[
\xi(t) := tv, \quad t \ge 0.
\]
The random set covered by time $t$ is 
\[
\mathcal W(t) = \bigcup_{x \in \Phi} ( x+\{vs,\, 0 \le s \le t\} + K).
\]
Letting $\xi^{(t)} := \{\xi(s),\, 0 \le s \le t\}$,
and $G^{(t)} :=  \xi^{(t)} + D$, we have
\begin{multline}
\label{CapIne}
\CAPRAND_{G^{(t)}}(K-x) 
= \P\big((\xi^{(t)} + D) \cap (K-x) \neq \varnothing\big)
\\
= \P\big(\{x+sv,\, 0 \le s \le t\}\,  \text{ intersects }K-D\big),
\end{multline}
which is laborious (but not impossible) to compute explicitly
(but we don't need the explicit formula).
Assume that the particles are initially placed at the points
of a Poisson process with intensity measure $\lambda(\d x)$.
Then, from \eqref{byfubini}, the distribution of
the detection time of $K$ is
\[
\P(S_K > t) =
\exp\bigg(-\E\, \lambda\bigg( \bigcup_{0 \le s \le t} [K-D-sv]\bigg)\bigg).
\]
Take now $K= r \B$, $D=r_0 \B$, with $\B$ the unit ball in $\R^d$.
Then $K-D = (r+r_0) \B$. Let 
\[
R:=r+r_0.
\]
Assuming further that $\lambda$ is isotropic
(invariant under rotations), so that
\[
\lambda(\d x) = \d \theta\, |x|^{d-1}\, \mu(\d |x|),
\]
where $\d\theta$ is the natural spherical measure on the boundary of the 
unit ball and $\mu$ a measure on $\R_+$,
we have
\[
\lambda\bigg( \bigcup_{0 \le s \le t} (K-D-sv) \bigg)
= \frac{1}{2} \lambda(R \B) + \frac{1}{2} \lambda(R \B + t v e_1)
+ \lambda(C_{R,tv}),
\]
where $C_{R,tv}$ is a cylinder with height $tv$. Specifically,
$C_{R,tv}$ contains the $x \in \R^d$ such that 
\begin{align*}
&x \cdot v \le t |v|^2, 
\\
&|v|^2 x - (x \cdot v) v \in |v|^2 \in R \B.
\end{align*}
Even when the law of $v$ is assumed isotropic,
the integrals can be quite hard to compute exactly unless we
further assume invariance under translations for $\lambda$, i.e., take
now $\lambda$ to be a multiple of the Lebesgue measure:
\[
\lambda(\d x) = \lambda \cdot \d x.
\]
In this case, things are very simple:
\begin{equation}
\label{verysimple}
\lambda\bigg( \bigcup_{0 \le s \le t} (K-D-sv) \bigg)
= \lambda \cdot R^d \omega_d + \lambda \cdot R^{d-1} \sigma_{d-1}
t |v|.
\end{equation}
Therefore, if $K-D = R \B$, and if the original location of
particles is a homogeneous Poisson process with intensity $\lambda$
then, assuming that $\E|v|<\infty$,
\begin{equation}
\label{lawinertial}
\P(S_K > t) = \exp(-\lambda R^d \omega_d) \,
\exp(-\lambda R^{d-1} \sigma_{d-1} (\E |v|) t).
\end{equation}
This gives
\begin{equation}
\label{expectationinertial}
\E S_K = \frac{e^{-\lambda \omega_d R^d}}
{\lambda d \omega_d (\E|v|) R^{d-1}}
= \frac{1}{\lambda d \omega_d (\E|v|)} ~ \frac{1}{R^{d-1}}
- \frac{1}{d} R + o(R), \quad \text{as } R \downarrow 0.
\end{equation}
However, if $\E|v|=\infty$, then the Poisson process \eqref{thinnedpoisson}
has intensity measure proportional to \eqref{CapIne} 
which integrates, with respect to the $d$-dimensional Lebesgue measure,
to the expectation of \eqref{verysimple} and this is infinity if $t >0$.
Therefore, 
\[
\text{ if  $\E|v|=\infty$ then $S_K=0$, a.s.}
\]
A more elaborate problem is the computation of the law of
\[
\inf\{t \ge 0: r_0 \B \cap  \mathcal W(t) \neq \varnothing,~
a \not \in \mathcal W(t)\}
\]
where
\[
\mathcal W(t) = \bigcup_{x \in \Phi} ( x+\{vs,\, 0 \le s \le t\} + r \B),
\]
which is the  set covered up to time $t$ by the inertial Boolean model
when the particles carry balls of radii $r$ each.
In other words, the problem is that of finding information about the first time
that the particles will detect a fixed ball of radius $r_0$
before anyone of them goes close to some point $a$ (the enemy).
This problem will be addressed in the future.

\section{The Brownian Boolean model}

We now specialize further and take $\xi$ to be a standard
Brownian motion in $\R^d$. In other words, $\xi(0)=0$
and $\xi=(\xi_1, \ldots, \xi_d)$, where the $\xi_i$, $i=1, \ldots, d$
are i.i.d.\ standard Brownian motions in $\R$.
Let $D$ and $K$ be balls with radii $r$ and $r_0$, respectively,
as in \eqref{Kset}, and let $S$ be the detection time of $K$.
The distribution of $S$ depends on $r_0$ and $r$ through their sum, so we write
\[
R=r_0+r
\]
for convenience.
Let
\[
\varrho(t) :=|\rho e_1 +\xi(t)|.
\]
Then $\varrho$ is standard Bessel process
of dimension $d$ started at $\varrho(0)=\rho$,
denoted as BES$^d(\rho)$ by Revuz and Yor \cite{RY91}. 
It is a strong Markov process (in fact, a Feller diffusion)
satisfying the It\^o equation
\[
\varrho(t) 
= \varrho(0) + \beta(t) + \frac{d-1}{2} \int_0^t \frac{1}{\varrho(s)} \,\d s,
\]
where $\beta$ is a standard Brownian motion in $\R$.
Letting $T^\rho_R$ be the first hitting time of the closed ball
of radius $R$ centered at the origin by the Bessel process started at $\rho$,
we have, from \eqref{SKrot},
\begin{equation}
\label{StailWiener}
\P(S>t) = \exp\bigg(-\lambda \sigma_{d-1} \int_0^\infty \rho^{d-1}\,
\P (T^\rho_R\le t)\, \d \rho \bigg) =: e^{-\lambda V^R(t)}.
\end{equation}

The infinitesimal generator $\mathcal A$ of $\varrho$ is the radial
part of $\frac{1}{2} \Delta$, where $\Delta$ is the Laplacian on $\R^d$:
\[
\mathcal A f(\rho) =
\frac{1}{2} \frac{1}{\rho^{d-1}} \frac{\partial}{\partial \rho} \bigg(\rho^{d-1}
\frac{\partial f}{\partial \rho}\bigg)
= \frac{1}{2} f''(\rho) + \frac{d-1}{2\rho} f'(\rho),
\]
acting on $C^2$ functions.
Therefore, for $s>0$, the function
\[
L(\rho) := \E [e^{-s T^\rho_R} ]
\]
satisfies 
\[
\mathcal A L = sLa
\]
i.e., the ODE
\begin{equation}
\label{BesselODE}
\frac{1}{2} L''(\rho) + \frac{d-1}{2\rho} L'(\rho) = sL(\rho),
\quad R < \rho < \infty,
\end{equation}
with boundary conditions
\begin{equation}
\label{bc}
L(R+)=1,
\quad
\lim_{\rho \to \infty} L(\rho)=0.
\end{equation}
Change variables using
\begin{equation}
\label{COV}
L(\rho) = \rho^b \widetilde L(a \rho),
\end{equation}
for appropriate constants $a >0, b \in \R$. 
The ODE reduces to
\[
a^2 \rho^2 \widetilde L''(a\rho) + (d-1+2b) a\rho \widetilde L'(a\rho)
+(b^2-2b+bd -2s\rho^2) \widetilde L(a\rho)=0.
\]
Choosing
\[
b= \frac{d}{2}-1
\]
gives
\[
a^2 \rho^2 \widetilde L''(a\rho) +  a\rho \widetilde L'(a\rho)
-(b^2+2s\rho^2) \widetilde L(a\rho)=0.
\]
Letting 
\begin{equation}
\label{a2s}
a=\sqrt{2s}
\end{equation}
(and letting $x:=a\rho$), we obtain the following ODE
\begin{equation}
\label{modbesode}
x^2 \widetilde L''(x) + x \widetilde L'(x) - (b^2+x^2) \widetilde L(x)=0.
\end{equation}
We recognize this as the modified Bessel ODE \cite[Sec.\ 3.7]{JACK}
the fundamental solutions of which are the modified
Bessel functions $I_{\pm b}$ and $K_{b}$.
The standard Bessel ODE differs from \eqref{modbesode}
by a change of sign in the last term. The fundamental solutions of the
standard Bessel ODE
are the Bessel functions of first and second kind $J_{\pm b}$, $N_b$,
whose series representations are easily obtained from the ODE; 
see equations (3.82), (3.83) and (3.85) in \cite{JACK}.
The modified Bessel functions $I_{\pm b}$ and $K_{b}$
(of first and second kind, respectively)
are related to $I_{\pm b}$ and $K_{b}$ via
\begin{align*}
I_{\pm b}(x) &= i^{-b} J_{\pm b}(ix)
\\
K_b(x) &= \frac{\pi}{2} i^b
[ i J_b(ix)- N_b(ix)],
\end{align*}
and are real-valued, despite appearances; see (3.100), (3.101)
and (3.86) in \cite{JACK}.
Since $I_{\pm b}$ explodes as $x\to\infty$, we are left with
only one choice for \eqref{modbesode}:
\[
\widetilde L(x) = C K_{b}(x),
\]
where $C$ is a constant.
In terms of the original function, i.e., using the change of
variables \eqref{a2s}, \eqref{COV},
\[
L(\rho) = C \rho^{-b} K(\rho\sqrt{2s}).
\]
The boundary conditions \eqref{bc} determine $C$:
\[
C= R^{b}/K_{b}(R \sqrt{2s}).
\]
So the solution to the ODE \eqref{BesselODE} 
with boundary conditions \eqref{bc}
is given by
\begin{equation}
\label{TLaplace}
L(\rho) = \E[e^{-sT^\rho_R} ]
= \frac{\rho^{-b} K_b(\rho\sqrt{2s})}
{R^{-b} K_b(R\sqrt{2s})}, \quad \rho \ge R.
\end{equation}

Compare now \eqref{StailWiener} with \eqref{detHOMmbm}. 
Since $\bigcup_{0 \le s \le t} (\xi(s) + D-K) = 
\bigcup_{0 \le s \le t} (\xi(s) + R \mathbb B)$, the quantity
in the  exponent
in \eqref{StailWiener} is the expected volume of the
Wiener sausage
\[
W^R(t) := \bigcup_{0 \le s \le t} (\xi(s) + R \mathbb B);
\]
see comments at the end of \S \ref{dtimbm}.
We have
\[
V^R(t):= \E \vol  W^R(t) = \int_{\R^d} \P(T^x_R \le t) \,\d x
= \omega_d R^d + 
\sigma_{d-1} \int_R^\infty \rho^{d-1} \P(T^\rho_R \le t)\,\d \rho.
\]
Via \eqref{TLaplace}, we have an expression for
the Laplace transform of the the expected volume of
the Wiener sausage:
\begin{align}
\widehat V^R(s) := \int_0^\infty e^{-st} V^R(t) \,\d t 
&=
\frac{\omega_d R^d}{s}
+ \frac{\sigma_{d-1}}{s} \int_R^\infty \big(\E e^{-s T^\rho_R}\big)
\rho^{d-1} \,\d \rho
\nonumber
\\
&= \frac{\omega_d R^d}{s}
+ \frac{\sigma_{d-1}}{s} \int_R^\infty
\frac{R^{\frac{d}{2}-1}}{\rho^{\frac{d}{2}-1}}
\frac{K_{\frac{d}{2}-1}(\rho\sqrt{2s})}{K_{\frac{d}{2}-1}(R\sqrt{2s})}
\rho^{d-1} \,\d \rho
\nonumber
\\
&= \frac{\omega_d R^d}{s}
+ \frac{\sigma_{d-1} R^{\frac{d}{2}-1}}{sK_{\frac{d}{2}-1}(R\sqrt{2s})}
\int_R^\infty K_{\frac{d}{2}-1}(\rho\sqrt{2s}) \rho^{d/2} \,\d \rho.
\label{VolLap}
\end{align}
We now use a couple of facts about the functions $K_b$; see \cite{AS}.
First, we have the recursion formula
\[
K_{b+1}(x) - K_{b-1}(x) = \frac{2b}{x} K_b(x).
\]
Second, we have the derivative
\[
K_b'(x) = \frac{b}{x} K_b(x) - K_{b+1}(x).
\]
Combining these we get
\[
\frac{\d }{\d x} \big( K_b(x) x^b \big) = - K_{b-1}(x) x^b,
\]
and so
\[
\int_x^\infty K_{\frac{d}{2}-1}(y) y^{\frac{d}{2}} \,\d y 
= K_{\frac{d}{2}}(x) x^{\frac{d}{2}},
\]
and, for $\lambda > 0$,
\[
\int_x^\infty K_{\frac{d}{2}-1}(\lambda y) y^{\frac{d}{2}} \,\d y 
=\frac{1}{\lambda} K_{\frac{d}{2}}(\lambda x) x^{\frac{d}{2}}.
\]
The last integral in \eqref{VolLap} evaluates to
\[
\frac{R^{\frac{d}{2}}}{\sqrt{2s}} K_{\frac{d}{2}}(R\sqrt{2s}),
\]
and so
\[
\widehat V^R_d(s) =  \frac{\omega_d R^d}{s}
+ \frac{\sigma_{d-1} R^{d-1}}{\sqrt{2s^3}} 
\frac{K_{\frac{d}{2}}(R\sqrt{2s})}{K_{\frac{d}{2}-1}(R\sqrt{2s})},
\]
where we added a subscript $d$ to indicate dependence on the dimension.
We can save some space by observing that, due to Brownian scaling,
\begin{equation}
\label{Vscaling}
V^R_d(t) = R^d V^1_d(t/R^2),
\end{equation}
\[
\widehat V^R_d(s) = R^{d+2} \widehat V^1_d(R^2 s).
\]
and so it is only
\begin{equation}
\label{itisonly}
\widehat V^1_d(s) =  \frac{\omega_d }{s}
+ \frac{\sigma_{d-1} }{\sqrt{2s^3}} 
\frac{K_{\frac{d}{2}}(\sqrt{2s})}{K_{\frac{d}{2}-1}(\sqrt{2s})}
\end{equation}
we should be looking for.
Since $K_b=K_{-b}$, the case $d=1$ is trivial:
\begin{equation}
\label{V11}
\widehat V^1_1(s) = \frac{2}{s} + \frac{2}{\sqrt{2s^3}}.
\end{equation}
Inverting this Laplace transform gives
\[
V_1^1(t) = 2 + \sqrt{\frac{8 t}{\pi}}.
\]
By the  scaling relation \eqref{Vscaling},
\[
V_1^R(t) = 2R + \sqrt{\frac{8 t}{\pi}},
\]
that is, the expected change of volume (=length) from its
initial value does not depend on $R$.

\subsection{Odd dimensions}
\label{odddimensions}
Consider now the case where
\[
d=2n+1, \quad n =0,1,\ldots
\]
We will produce an algorithm for computing $\widehat V^1_{2n+1}(s)$ recursively,
and carry out its first few steps.
The modified Bessel functions of half-integer order have a simple
form:
\[
K_{n+\frac{1}{2}}(x) = \sqrt{\frac{\pi}{2}} \frac{e^{-x}}{\sqrt{x}}
y_{n}(1/x), \quad n=0,1,\ldots,
\]
where $y_n(x)$ is a polynomial of degree $n$ with integer coefficients:
\[
y_n(x) := \sum_{k=0}^n \frac{(n+k)!}{(n-k)! k!} \bigg(\frac{x}{2}\bigg)^k,
\quad n=0,1,\ldots,
\]
known as the Bessel polynomial of degree $n$; see \cite[formula (3)]{KF}.
Note, in particular, that
\[
y_n(x) = (2n-1)!! x^n + (2n-1)!! x^{n-1} + \cdots + \frac{n(n+1)}{2} x + 1,
\]
i.e., the coefficients of the two highest powers are equal to the double factorial
\[
(2n-1)!! = \frac{(2n)!}{n! 2^n} = (2n-1)(2n-3)(2n-5) \cdots 3 \cdot 1.
\]
See Appendix \ref{modbes} for a table of the first few Bessel polynomials
and their corresponding Bessel functions.
Consequently,
\begin{equation}
\label{Vhat}
\widehat V^1_{2n+1}(s) = \frac{\omega_{2n+1}}{s} 
+\frac{\omega_{2n+1}}{s}~ \frac{2n+1}{\sqrt{2s}}~
\frac{y_n(1/\sqrt{2s})}{y_{n-1}(1/\sqrt{2s})},
\quad n=0,1,\ldots,
\end{equation}
where
\[
y_{-1}(x) :=1,
\]
as follows by comparison to \eqref{V11}. We thus have
\begin{align*}
\frac{1}{\omega_3} \widehat V^1_3(s)
&= \frac{1}{s} +\frac{3}{s} ~\frac{\sqrt{2s}+1}{2s}
\\
\frac{1}{\omega_5} \widehat V^1_5(s)
&= \frac{1}{s} + \frac{5}{s}~ \frac{(2s)^{3/2} + 6s + 3 \sqrt{2s}}
{(2s)^{3/2}~[\sqrt{2s}+1]}
\\
\frac{1}{\omega_7} \widehat V^1_7(s)
&= \frac{1}{s} + \frac{7}{s}~ 
\frac{(2s)^{3/2} + 12 s + 15 \sqrt{2s} + 15}
{\sqrt{2s}~ [(2s)^{3/2} + 6s + 3 \sqrt{2s}]}
\\
\frac{1}{\omega_9} \widehat V^1_9(s)
&= \frac{1}{s} + \frac{9}{s}~ 
\frac{(2s)^{5/2} + 40 s^2 + 45 (2s)^{3/2} +210 s + 105\sqrt{2s}}
{(2s)^{3/2}~ [(2s)^{3/2} + 12 s + 15 \sqrt{2s} + 15]}.
\end{align*}
These Laplace transforms can, in principle, be inverted by
using partial fraction expansion and the fact that (see
Erd\'elyi {\em et al.} \cite[Ch.\ 7, p.\ 233]{ERD})
\[
\frac{1}{\sqrt{s}+\beta} =
\int_0^\infty e^{-st} 
\bigg[\frac{1}{\sqrt{\pi t}} - \beta e^{\beta^2 t} \erfc(\beta \sqrt{t})\bigg]
\,\d t,
\]
where
\[
\erfc(t):=1-\erf(t),\quad
\erf(t) := \frac{2}{\sqrt{\pi}} \int_t^\infty e^{-u^2} \,\d u.
\]
For example, writing
\[
\frac{1}{\omega_3} \widehat V^1_3(s) 
= \frac{1}{s} + \frac{3}{\sqrt{2s^3}} + \frac{3}{2s^2},
\]
we obtain
\[
\frac{1}{\omega_3}V^1_3(t) = 1 + \frac{6}{\sqrt{\pi}}\sqrt{t} +\frac{3}{2} t.
\]
Expanding $\widehat V^1_5(s)$, we obtain
\[
\frac{1}{\omega_5} \widehat V^1_5(s)
=\frac{1}{s} + \frac{10}{2s}~
\bigg[\frac{1}{\sqrt{2s}+1} + \frac{3}{\sqrt{2s}(\sqrt{2s}+1)}
+ \frac{3}{2s(\sqrt{2s}+1)} \bigg].
\]
Since
\begin{align*}
\frac{1}{s(\sqrt{s}+1)} &= \int_0^\infty e^{-st} [1-e^t \erfc(\sqrt{t})] \,\d t
\\
\frac{1}{s\sqrt{s}(\sqrt{s}+1)} = \frac{1}{s\sqrt{s}} - \frac{1}{s(\sqrt{s}+1)}
&=
\int_0^\infty e^{-st}
\bigg[2\sqrt{\frac{t}{\pi}} -1 + e^t \erfc(\sqrt{t}) \bigg] \,\d t
\\
\frac{1}{s^2(\sqrt{s}+1)} &=
\int_0^\infty e^{-st} \bigg[
1+t - 2\sqrt{\frac{t}{\pi}} - e^t \erfc(\sqrt{t}) 
\bigg] \,\d t,
\end{align*}
letting $g_1(t), g_2(t), g_3(t)$ be the functions in the square brackets
of the last three lines,
we have
\begin{align*}
\frac{1}{\omega_5} V^1_5(t)
= 1 + 10 
\bigg[\frac{1}{2} g_1(t/2) + \frac{3}{2} g_2(t/2) + \frac{3}{2} g_3(t/2)\bigg]
&=
6 -5 e^{t/2} \erfc(\sqrt{t/2})+ \frac{15}{2} t. 
\end{align*}
Using the scaling relation \eqref{Vscaling}, we can now
obtain $V^R_d(t)$, for $d=1,3,5$ and, therefore,
the distribution of the detection time via \eqref{StailWiener}.
\begin{align}
d=1:&  \qquad \P(S> t) = \exp\big(-2\lambda R -4\lambda\sqrt{t/\pi}\big).
\label{d1}
\\
d=3: & \qquad
\P(S>t) = \exp\big(-\frac{4\pi\lambda}{3}R^3
-8\sqrt{\pi t}\lambda R^2 - 2\pi \lambda t R\big).
\label{d2}
\\
d=5: & \qquad
\P(S>t) = \exp\big(-\frac{16 \pi^2}{5}R^2 
+ \frac{8R^5}{3} e^{t/2R^2} \erfc(R^{-1} \sqrt{t/2}) 
-4\pi^2 R^3 t\big) 
\label{d3}
\end{align}
Notice that the exponent is not a polynomial in $R$,
as is apparent for the $d=5$ case.\footnote{Answering a question posed
by G\"unter Last to me a few years ago.}

A more efficient way of doing the above is by means of the recursion formula
\[
y_n(x) = (2n-1) x y_{n-1}(x) + y_{n-2}(x), \quad n=1,2,\ldots,
\]
with initial conditions $y_{-1}(x)=y_0(x)=1$. See \cite[\S 7]{KF}.
Let 
\[
H_n(s) := \frac{\sqrt{2s}}{2n+1} \bigg(s \frac{\widehat V^1_{2n+1}(s)}{\omega_{2n+1}}
-1 \bigg).
\]
From \eqref{Vhat} we have
\[
H_n(s) 
= \frac{y_n(1/\sqrt{2s})}{y_{n-1}(1/\sqrt{2s})},
\]
and so, from the recursion formula for Bessel polynomials,
\[
H_n(s) = \frac{2n-1}{\sqrt{2s}} + \frac{1}{H_{n-1}(s)},
\quad n=1,2,\ldots,
\]
where $H_0(s)=1$.
This gives a recursive way for computing $\widehat V^1_{2n+1}(s)$.

Let
\[
a_n(s) := \frac{2n-1}{\sqrt{2s}}.
\]
A ``closed'' formula can also be obtained:
\begin{align}
\label{contfrac}
\frac{\widehat V^1_{2n+1}(s)}{\omega_{2n+1}}
&= \frac{1}{s} +  \frac{a_{n+1}(s)}{s}
\,
\left\{ 
a_n(s) + \cfrac{1}{a_{n-1}(s)
	+ \cfrac{1}{a_{n-2}(s) + \cfrac{1} 
{\ddots +\frac{1}{a_1(s)+1} }}}
\right\}
\\
\label{contfrac2}
&= \frac{1}{s} + \frac{2n+1}{2s^2}\,
\left\{ 
(2n-1) + \cfrac{\sqrt{2s}}{(2n-3)
        + \cfrac{\sqrt{2s}}{(2n-5) + \cfrac{\sqrt{2s}} 
{\ddots +\frac{\sqrt{2s}}{1+\sqrt{2s}} }}}
\right\}
\end{align}

\subsection{Large time asymptotics in all dimensions}
\label{largetime}
Expression \eqref{itisonly} allows us to find logarithmic asymptotics
for $\P(S > t)$, as $t \to \infty$, in any dimension $d$. We 
repeat the expression here:
\[
\frac{1}{\omega_d}\widehat V^1_d(s) =  \frac{1}{s}
+ \frac{d}{\sqrt{2s^3}} 
\frac{K_{\frac{d}{2}}(\sqrt{2s})}{K_{\frac{d}{2}-1}(\sqrt{2s})} =: 
\frac{1}{s} + \widehat g_d(s).
\]
The following asymptotics are known for Bessel functions \cite{AS}. For
any $b > 0$, as $z \to 0$,
\[
K_b(z) \sim 2^{b-1} \Gamma(b) z^{-b}.
\]
Therefore, for $d \ge 3$, we have
\[
\widehat g_d(s) \sim \frac{d(d-2)}{2s^2}, \quad s \to 0.
\]
Hence
\[
g_d(t) \sim \frac{d(d-2)}{2} t, \quad t \to \infty.
\]
By the scaling equation \eqref{Vscaling} and expression
\eqref{StailWiener}, we obtain
\begin{align*}
\P(S>t) &= \exp(-\lambda \omega_d R^d -\lambda \omega_d R^d g_d(t/R^2)).
\\
&\stackrel{\log}{\text{\Large $\sim$}} 
\exp\bigg(-\lambda \omega_d \frac{d(d-2)}{2} R^{d-2} t\bigg),
\quad d \ge 3.
\end{align*}
The case $d=2$ has to be treated differently as it requires the
behavior of $K_0$ near zero which is different:
\[
K_0(z) \sim \log(1/z) , \quad z \to 0.
\]
Hence
\[
\widehat g_2(s) = \frac{2}{\sqrt{2s^3}}
\frac{K_1(\sqrt{2s})}{K_0(\sqrt{2s})}
\sim
\frac{2}{\sqrt{2s^3}} 
\frac{1}{ \sqrt{2s }\log \frac{1}{\sqrt{2s}}}
=\frac{2}{s^2 \log(1/s)}
\]
Notice that 
\[
\widehat g_2(s) = s^{-2} \ell(1/s),
\]
where the function
$\ell(z) = 2/\log(z)$ is slowly varying at infinity, viz.,
$\ell(\kappa z)/\ell(z) \to 1$, as $z \to \infty$, for all $\kappa >0$.
Combining Karamata's Tauberian theorem with the monotone 
density theorem ($g_2(t)$ is increasing function of $t$) 
we conclude that 
\[
g_2(t) \sim  t \ell(t) = \frac{2t}{\log t},\quad \text{  as $t \to \infty$.}
\]
Arguing as before, this means that, as $t \to \infty$,
\begin{align*}
\P(S>t) &= \exp(-\lambda \omega_2 R^2 -\lambda \omega_2 R^2 g_2(t/R^2))
\\
&\stackrel{\log}{\text{\Large $\sim$}} 
\exp\bigg(-\frac{2\pi\lambda t}{\log t}\bigg).
\end{align*}
This agrees with the result of \cite[Theorem 2]{SPI}.

\subsection{Expectations}
For $d=1$, using $\int_0^\infty \exp(-\sqrt{t}) \,\d t=2$, we can compute 
the expectation of $S$ explicitly by integrating 
\eqref{d1}:
\[
\E S := \frac{\pi}{8} ~ \frac{e^{-2\lambda R}}{\lambda^2}.
\]

For $d=2$, we have no explicit expression, but the earlier asymptotic
expression for large $t$ tells us that, as $R \to 0$, $\E S$
converges to a constant. This is a manifestation of the fact
that Brownian motion is neighborhood recurrent in 2 dimensions.

For $d=3$, we can use the integral
$\int_0^\infty \exp(-\sqrt{t}-t) \,\d t=
1+ \frac{\pi e^{1/4}}{2}(\erf(1/2)-1)$ to integrate \eqref{d3}:
\[
\E S =
\frac{e^{-\frac{4}{3} \pi \lambda R^3}}{ 2(\pi\lambda R)^{3/2}}
\bigg[
\sqrt{\pi\lambda R} + 2\sqrt{2}\pi\lambda R^2 e^{8\lambda R^3}
\big( \erf( 2\sqrt{2\lambda} R^3/2 )-1 \big)
\bigg]
\]

For higher dimensions, we can translate the previous asymptotics
for $t \to \infty$ into asymptotics for $R\downarrow 0$ and obtain that
\begin{equation}
\label{expectationbrownian}
\E S \sim \frac{c_d}{R^{d-2}}, \text{ as } R \downarrow 0,
\end{equation}
where $c_d$ is a constant depending on $d$ only. This
estimate holds for all $d \ge 2$.
Comparing \eqref{expectationbrownian} with \eqref{expectationinertial}
we find that in any dimension $d\ge 2$, it is better to
make particles (sensors) perform Brownian motions
rather than random straight lines with finite mean velocity
if the goal is to detect a small object.\footnote{This answers
a question posed by Venkat Anantharam a few years ago, who, jokingly,
commented that an engineer would never have sensors perform 
Brownian motions.}

%
%

\section{Concluding remarks and open problems}
In this paper, we reviewed the mobile Boolean model, focusing,
in particular, in the inertial and Brownian cases. 
For the inertial case, we have an explicit expression 
\eqref{lawinertial}
for the distribution of the detection
time of in any dimension.
For the Brownian case, we have explicit expressions
\eqref{d1}, \eqref{d2}, \eqref{d3} in dimensions $d=1$, $3$ and $5$,
an algorithm for computing an explicit expression when $d$ is odd
(Section \ref{odddimensions})
and logarithmic asymptotics when $d$ is even (Section \ref{largetime}).
We worked with a target set $K$ which is a ball. This enabled us to reduce
the a multi-dimensional problem to one dimension. 
We also gave formulas or estimates for expectations.

For general compact set $K$,
and dimension $d=3$,
Spitzer's \cite{SPI} paper gives asymptotic estimates for
the expected volume of a $K$-Wiener sausage,
in terms of the Newtonian capacity $\operatorname{cap}(K)$ of $K$,
using entirely probabilistic methods:
\[
\E \vol(\xi^{(t)} + K) =
\vol(K) + \operatorname{cap}(K)\,t 
+ 4(2\pi)^{-3/2} \operatorname{cap}(K)^2\, \sqrt{t}
+o(\sqrt{t}), \quad d=3.
\]
Translating this into a detection time probability estimate,
and using the scaling relation, we have
\[
P(S_{R K} > t)
= \exp\big(
\vol(K)\,R^3 +  \operatorname{cap}(K)\,R^2 t
+  4(2\pi)^{-3/2} \operatorname{cap}(K)^2\, R \sqrt{t}
+o(R \sqrt{t})
\big).
\]

We next pose some open problems.

\paragraph{Open problem 1.}
We did not touch at all the coverage problem, i.e., the law of
the random variable $\inf\{t \ge 0:~ K \subset \mathcal W(t)\}$.
For the Brownian Boolean model, and when the diameter of $R$ tends to infinity,
the problem has been solved in \cite{PSSS}. What are the corresponding 
asymptotics for the inertial cases?

\paragraph{Open problem 2.}
Let $K_1$ and $K_2$ be two sets (e.g., balls of radii $r_1, r_2$),
and let $S_{K_1}, S_{K_2}$ be their detection times.
Find the probability $P(S_{K_1} < S_{K_2})$.

\paragraph{Open problem 3.}
For the inertial Boolean model, find the distribution of the
detection time of a ball before a fixed point is hit.
(See remarks at the end of Section \ref{inbm}.)

\paragraph{Open problem 4.}
Investigate further the algorithm of Section \ref{odddimensions}
and, in particular, the ``closed'' formula \eqref{contfrac}-\eqref{contfrac2}.

\paragraph{Open problem 5.}
Let there be an independent space-time
Poisson process $\Psi$ in $\R^d \times \R_+$ with fixed intensity.
Interpret its points as ``customers''. The Brownian Boolean model
is a space-time serving mechanism clearing points whenever it meets them.
Find necessary and sufficient stability conditions.
This problem is related to a number of recent stability problems in
queueing theory where the spatial dimension is just as important as
the time dimension \cite{SF,BFL,RFK}.

\appendix
\section{Modified Bessel functions of second kind of
half-integer order and their
corresponding Bessel polynomials}
\label{modbes}
\begin{align*}
K_{n+\frac{1}{2}}(x)  &= \sqrt{\frac{\pi}{2}} \frac{e^{-x}}{\sqrt{x}}
& y_n(x)   &=   \sum_{k=0}^n \frac{(n+k)!}{(n-k)! k!} \bigg(\frac{x}{2}\bigg)^k
\\
& & \\
K_{1/2}(x) &= \sqrt{\frac{\pi}{2}} \frac{e^{-x}}{\sqrt{x}}
&y_0(x) &= 1 \\
K_{3/2}(x) &= \sqrt{\frac{\pi}{2}} \frac{e^{-x}}{x^{3/2}} (x+1) 
&y_1(x) &= 1+x \\
K_{5/2}(x) &= \sqrt{\frac{\pi}{2}} \frac{e^{-x}}{x^{5/2}} (x^2+3x+3)
&y_2(x) &= 1+3x+3x^2 \\
K_{7/2}(x) &= \sqrt{\frac{\pi}{2}} \frac{e^{-x}}{x^{7/2}} (x^3+6x^2+15x+15)
&y_3(x) &= 1+6x+15x^2+15x^3\\
K_{9/2}(x) &= \sqrt{\frac{\pi}{2}} \frac{e^{-x}}{x^{9/2}} 
(x^4+10x^3+45x^2+105 x+105)
&y_4(x)&=1+10 x+ 45 x^2 + 105 x^3 + 105 x^4.
\end{align*}

\vspace*{1cm}
\hfill
\noindent
\begin{minipage}[t]{6cm}
\small \sc
Takis Konstantopoulos
\\
Department of Mathematics
\\
Uppsala University
\\
751 06 Uppsala
\\
Sweden
\\
{\tt takiskonst@gmail.com
\\
www.math.uu.se/$\sim$takis}
\end{minipage}

\end{document}